\documentclass[a4paper,12pt]{article}
\usepackage{amsmath}
\usepackage{amsthm}
\usepackage{amssymb}
\usepackage[authoryear]{natbib}
\usepackage[T1]{fontenc}
\usepackage[english]{babel}
\usepackage{amsfonts}
\usepackage{graphicx,color,float}

\usepackage[latin1]{inputenc}  

\usepackage[colorlinks=true]{hyperref}

\title{\bf Decision-making and interacting neuron populations}
\author{S. MANCINI\\
MAPMO, Université d'Orléans}
\date{}

\begin{document}

\maketitle 
\begin{abstract}In this article we present the modeling of bi-stability view problems described by the  activity or firing rates of two interacting population of neurons. Starting from the study of a complex system, the system of stochastic differential equations describing the time evolution of the activity of the two populations of neurons, we point out the strength and weakness of this model and consider its associated partial differential equation, which resolution gives statistical information on the firing rates distributions. The slow-fast characterization of the solutions finally leads us to a complexity reduction of the model by the definition of a one-dimensional stochastic differential equation and its associated one-dimensional partial differential equation. This last model turns out to be well adapted to the resolution of the problem giving access, in particular, to reaction times and performance, two macroscopic variables describing the decision-making in the view problem.
\end{abstract}

\section{Introduction}
Decision-making problems in social and natural sciences are often described by means of complex systems governed by differential equations giving the time variation of some  quantities. We can represent the two choices decision making situation like a set of particles evolving in double wells potentials (potential functions with two minima or stable equilibrium points and one maximum or unstable equilibrium point) and submitted to interactions. Each well represents one of the decision states and corresponds to one attractor of the system. The function describing the double well potential is usually a fourth order polynomial (like in the Van der Pool equation), and the problem can be explicitly solved, see for example the book by \citet{Galam} for an application to social sciences. The decision-making process may also be described by the evolution of the reaction times and performance, two macroscopic variables representing the mean minimal time a subject needs to make a decision (or a particle needs to exit a potantial well), and the amount of subjects having choosen a particular decision state  (or the sum of all particles being in a potential well) at a given time.

In this paper we deal with bi-stability visual situations. The decision-making process in this context involve a huge number of interacting neurons and it is not possible to describe it by the knowledge of each single neuron. The synchronization of the neurons activity leads to an equilibrium representing the decision. We can briefly sketch the situation as follows.  A subject is asked to choose between two possible views of a picture. His sight has to focus on one of these views, and this is done at a neuronal level. The decision is taken once the focus is done. Neuron-physicists are then interested in the two macroscopic quantities: the reaction times and the performance. Since neurons in the visual cortex have different skills and are connected, this problem can't be modeled by the description of a single neuron activity, but different populations of neurons in interaction must be considered. 
In computational neurosciences, the decision making of interacting population of neurons (excitatory and inhibitory ones) have been successfully described by a system of deterministic differential equations, called the Wilson-Cowan system, see the seminal work by \citet{WC}. In this model, the unknowns are function of time only and represent the mean firing rate of each population of neurons, i.e. for each population, the mean frequency of the neuronal signal, hence its activity. Moreover, the underlying potential is not a fourth order polynomial function and can't be explicitly computed. More recently, noise has been added to the model, see \citet{DM}, in order to account for  the finite number of neurons in the mean field approximation used to derive the Wilson-Cowan model. The non-linearity in the model makes its mathematical analysis difficult. In particular, it is not possible to write the explicit solution of the stationary associated problem. Nevertheless,  in \citet{DM} the authors numerically show, applying the moment analysis, that for the ranges of  parameters they are interested in, the solutions are bi-modal, i.e. double peaked. This method works well, but no closure to the system of equation is provided. In order to write an approximation of the explicit stationary solution to the problem, knowing that the solutions must be bi-modal and applying Taylor expansion methods, it is also possible to define a fourth order polynomial $V$, passing from the equilibrium points, see for example \citet{RL}. This function $V$ is then also used to compute reaction times (by means of Kramers formula) and performance (by means of the steady state). This approach, which is usually applied in computational neurosciences, gives results in agreement with experimental data, see \citet{RL}, but only holds locally. 

In this paper, we will consider the partial differential equation associated to the stochastic system which describes the evolution of the probability distribution function in terms of the firing rates of the two neuron populations. We will see in the sequel how its mathematical analysis and numerical simulations can help to reduce the complexity of the problem leading to faster computations of the reaction times and the performance, and compare the results of the simplified model with the initial one. The presented complexity reduction method is a good way to overcome the difficulty of not knowing the explicit form of the underlying potential and the approximated solutions are defined on the whole domain we are interested in, and not only locally as usually done in computational neurosciences. The present work resumes the results of several papers done in collaborations with  J. A. Carrillo (London), G. Deco (Barcelona) and S. Cordier (Orléans), see \citet{CCM}, \citet{CCM2} and \citet{CCDM}. 

\section{The mathematical model}
Recently bi-stability visual problems have been investigated by considering systems of stochastic differential equations which describe the time evolution of the firing rates for two or more interacting populations of neurons (see for example \citet{DM} and \citet{RL}). This kind of models, based on the deterministic Wilson-Cowan one (see \citet{WC}), permits to numerically evaluate the subject reaction times and the performance together with their variations with respect to the differences on the applied stimuli and/or the weight of the interactions. For instance, reaction times correspond to the time needed for the subject to make a decision, and performance is the number of good responses taken by the subject without limitation on time.
The model can be interpreted from a physical point of view as particles trapped in a double (or multiple) well potential, reaction times corresponding then to the exit time from a well and performance being given by the density contained in the well associated to the correct answer.

The model studied in \citet{DM} considers the time evolution of the firing rates $\nu_1=\nu_1(t)$ and $\nu_2=\nu_2(t)$ of two neuron populations. Their behavior satisfies the following system of stochastic differential equations:
\begin{equation}\label{dyn}
\begin{cases}
\textrm{d} \nu_1 = \psi_1(\nu_1,\nu_2) \textrm{dt} + d\xi  \\
\textrm{d} \nu_2 = \psi_2(\nu_1,\nu_2) \textrm{dt} + d\xi,
\end{cases}\end{equation}
where $d\xi$ is a white noise of standard deviation $\beta$ and $\psi_1$, $\psi_2$ are the dynamical part of the equations and model the neuronal activity. They are defined by:
$$\psi_1= -\nu_1 + \phi (\lambda_1 + w\ \nu_1+\hat{w}\ \nu_2),$$
$$\psi_2= -\nu_2 + \phi (\lambda_2 + \hat{w}\ \nu_1+w\ \nu_2),$$
with $\phi(z)$ the so-called response function to the mean excitation $z$, defined by the sigmoid:
$$\phi(z)=\frac{\nu_c}{1+\exp(-\alpha(z/\nu_c -1))},$$
with $\nu_c$ and $\alpha$ parameters that are fixed by biology and where the mean excitation $z$ is given by the sum of the applied stimuli ($\lambda_1$ or $\lambda_2$) and the internal activities given by a linear combination of the activity of each population weighted respectively by $w$ and $\hat{w}$ depending if we are considering the same population or not. Note that, the weights being symmetric, if the applied stimuli are the same (i.e. $\lambda_1=\lambda_2$), then the problem is symmetric and is defined as the unbiased case, whereas if one of the stimuli is larger than the other (say $\lambda_1 =\lambda_2 +\Delta \lambda$), the problem loses its symmetry and we define this situation as the biased case, with the bias given by $\Delta \lambda$.  In the following numerical results, with the exeption of those in figure \ref{rt-perf} which consider a slightly different potential, $\nu_c=20$, $\alpha=4$, $w=0.45$, $\hat w=1.23$, $\beta=0.3$, $\lambda_1=\lambda_2=15$ and $\Delta \lambda=0.01$.

It's well-know in literature that we can deduce a Fokker-Planck equation from system \eqref{dyn} applying Ito calculus or considering the forward Kolmogorov equation associated to \eqref{dyn}: for $(t,\nu_1,\nu_2)\in (0,+\infty)\times \Omega$,
\begin{equation}\label{FP}
\partial_t p + \nabla \cdot \left( F\ p - \frac{\beta^2}{2} \nabla p \right) =0 ,
\end{equation}
where $p=p(t,\nu_1,\nu_2)$ is the probability distribution function representing the probability that at time $t\geq 0$ the firing rates are in $(\nu_1,\nu_2) \in \Omega \subset \mathbb{R}_+^2$, and with $F=F(\nu_1,\nu_2)=(f(\nu_1,\nu_2),g(\nu_1,\nu_2))$ the drift term. The domain $\Omega$ being bounded (the square $[0,\nu_m]\times [0,\nu_m]$, with $\nu_m$ the maximal firing rate value for the neuron populations), we complete equation \eqref{FP} by the following Robin type (or no flux) boundary conditions: on $\partial \Omega$,
$$F\ p - \frac{\beta^2}{2} \nabla p=0 \ , $$
and we finally consider the normalized initial condition:
$$p(0,\nu_1,\nu_2)=p_0(\nu_1,\nu_2) \geq 0.$$
As proven in \citet{CCM}, under the assumption of incoming flux, the problem is well posed and there exists a unique steady state solution of the stationary Fokker-Planck equation. Nevertheless, as explained in \citet{CCM} there is no potential function $V=V(\nu_1,\nu_2)$ such that $F=-\nabla V$. This fact implies that it is not possible to write explicitly the steady state associated to \eqref{FP}.  Recall that the steady state and the potential $V$ are essential for the computing of reaction times and performance. 

\begin{figure}[H]
\begin{center}
\includegraphics[width=6cm,height=6cm]{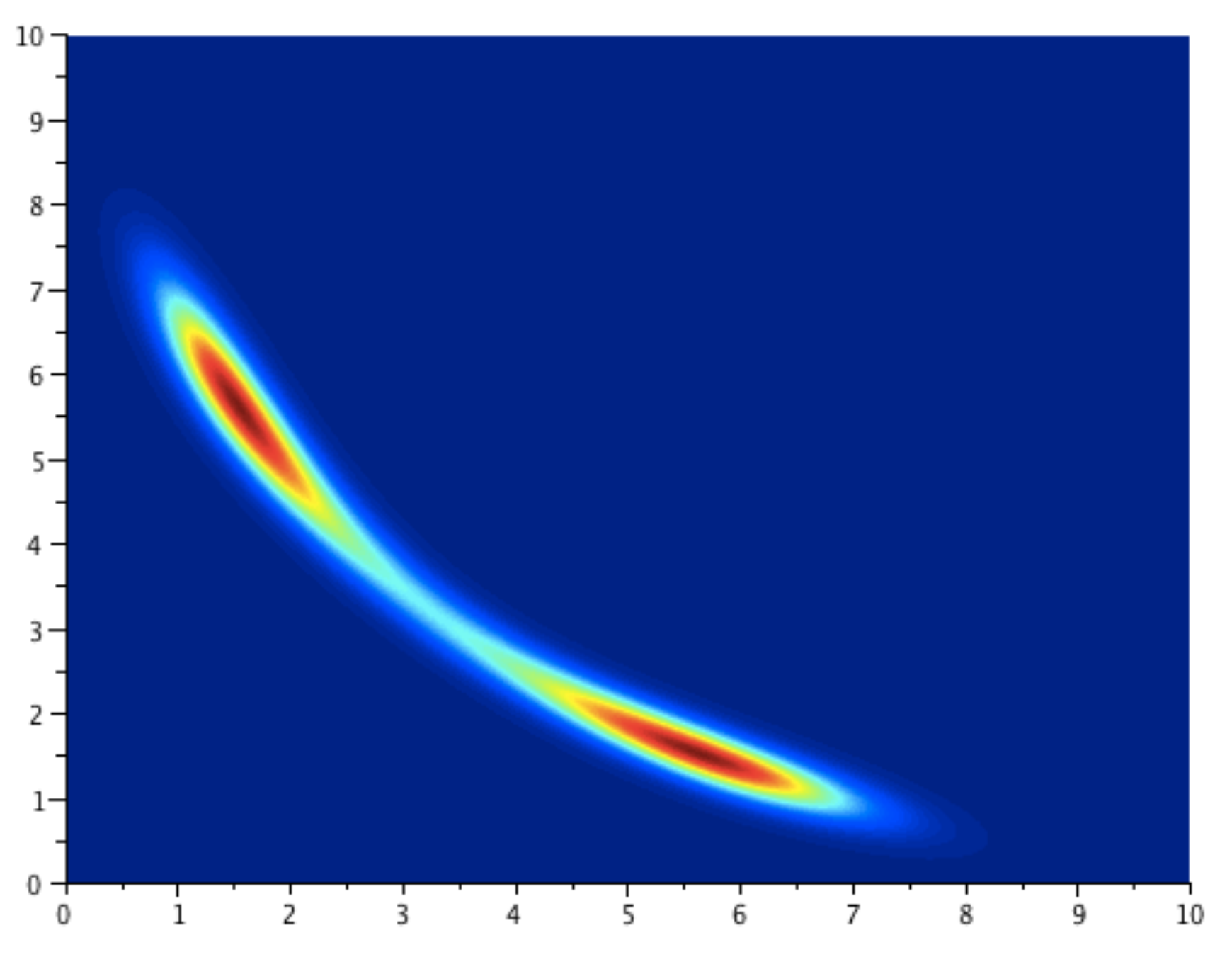}\; 
\includegraphics[width=6cm,height=6cm]{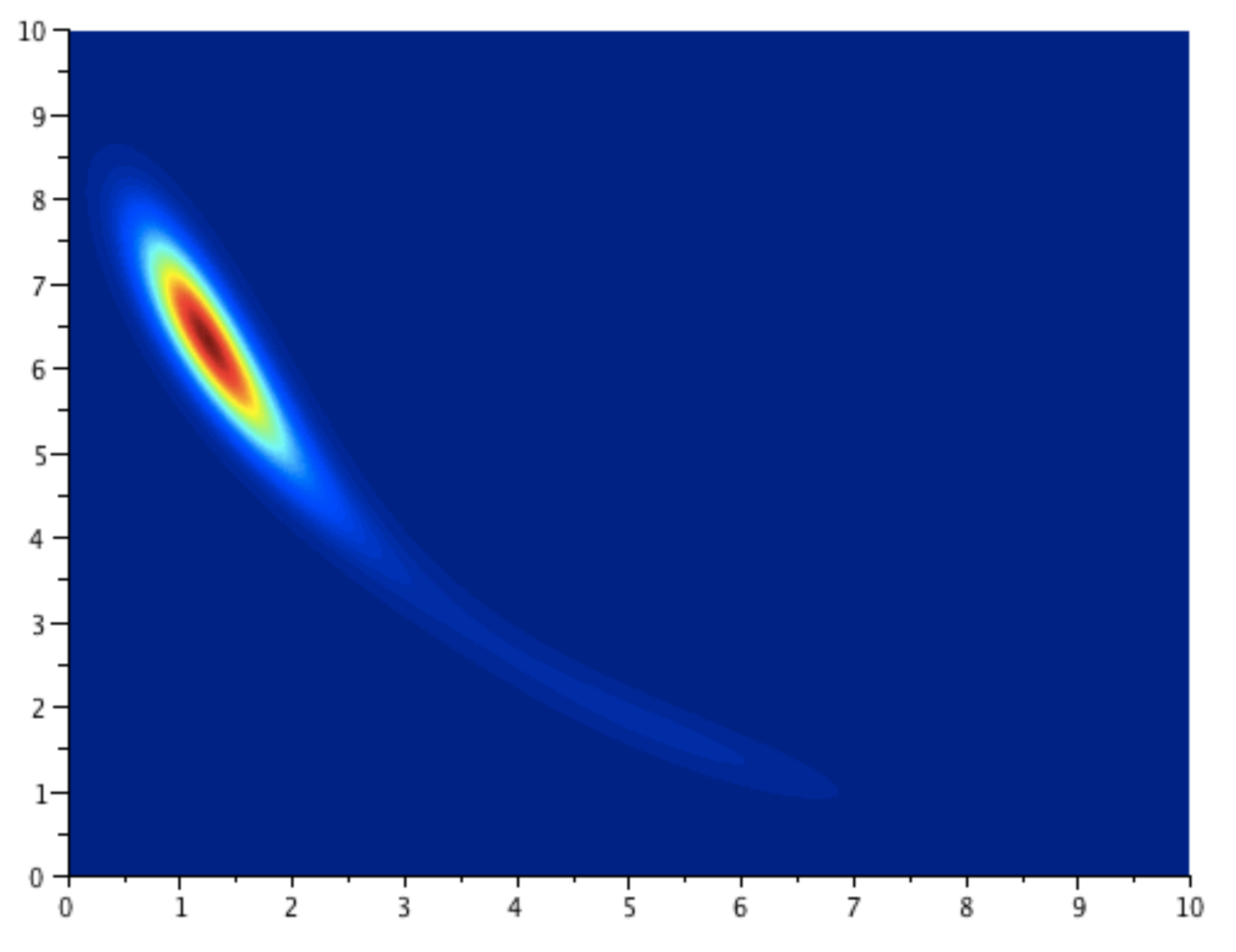}
\end{center}
\caption{Solution to equation \eqref{FP}. Left~: unbiased case. Right : biased case.}\label{contour}
\end{figure}

The bi-dimensional behavior of the solution of \eqref{FP} at a given time is shown in figure \ref{contour}: the solution is concentrated around the two stable equilibrium points and is aligned along the equilibrium manifold. The bi-modal aspect of the solution is well captured by the numerical simulations of the Fokker-Planck equation (left). When the situation is biased (i.e. one of the applied stimuli is bigger), then one of the wells is deeper than the other and the symmetry of the problem is lost. In this situation the solution at equilibrium (or for large times) is concentrated around the equilibrium point corresponding to the deeper well, (right). Note that in \citet{DM} the authors represented the solution by means of their marginals (i.e. the projections along each axis) and no bi-dimensional numerical result was obtained. 

\section{Complexity reduction}
Although the behavior of the solution of equation \eqref{FP} is in agreement with what is expected physically, the application to real problems of this equation, or the use of the associated numerical simulations, are not competitive, since nor the potential function $V$ nor the steady state are known and to approximate them numerically requires very long CPU times.

Nevertheless, the study of problem \eqref{FP} shows that its solution is characterized by a slow-fast  behavior: rapid diffusion towards the equilibrium manifold and slow drift along the manifold towards the stable equilibrium points. 
In \citet{CCM2} we have proposed a complexity reduction of \eqref{EDS} based on this slow-fast characteristic of the problem and leads to a one-dimensional Fokker-Planck equation living on the equilibrium manifold (see \citet{BG}). 
The fast convergence being along a direction which is given by a linear combination of $\nu_1$ and $\nu_2$, we can define two new variables $x$ and $y$ respectively as the variable along which the fast convergence is done and the one corresponding to the slow direction (see \citet{CCM2} for more details). With this change of variables the stochastic system \eqref{dyn} transforms in:
\begin{equation}\label{EDS}
\begin{cases}
\textrm{d}x = f(x,y) \textrm{dt} + d\xi_x ,\\
\textrm{d}y = g(x,y) \textrm{dt}+ d\xi_y ,
\end{cases}
\end{equation}
where $f(x,y)$, respectively $g(x,y)$, are the linear combination of the functions $\psi_1$ and $\psi_2$, and where $d\xi_x$ et $d\xi_y$ are two white noises of standard deviation $\beta_x$ and $\beta_y$. Summarizing, we may say that $f$ and $g$ are the functions describing the activity of the combined firing rates $x$ and $y$, respectively.

\begin{figure}[H]
\begin{center}
\includegraphics[width=.49\textwidth]{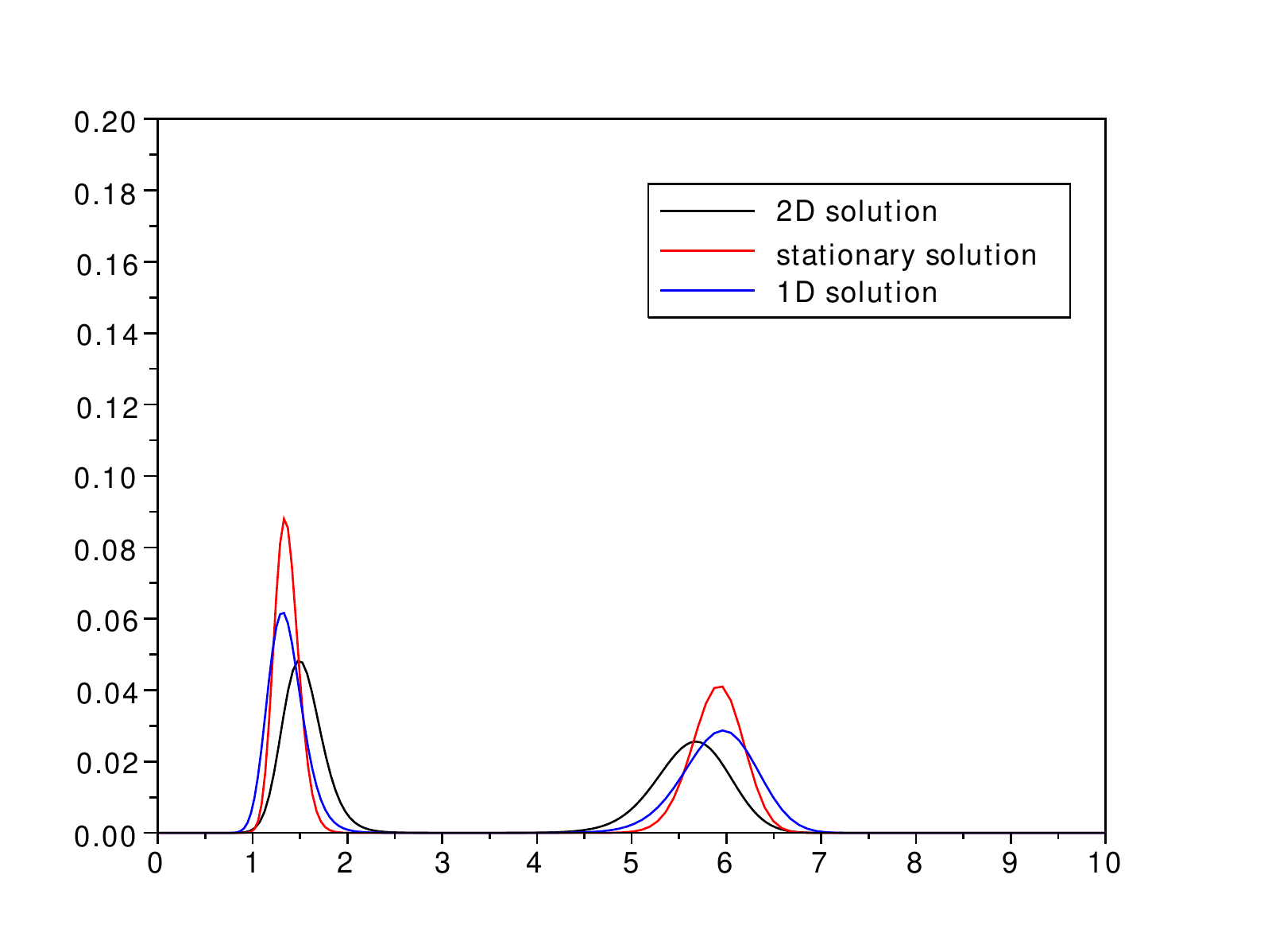}\; 
\includegraphics[width=.49\textwidth]{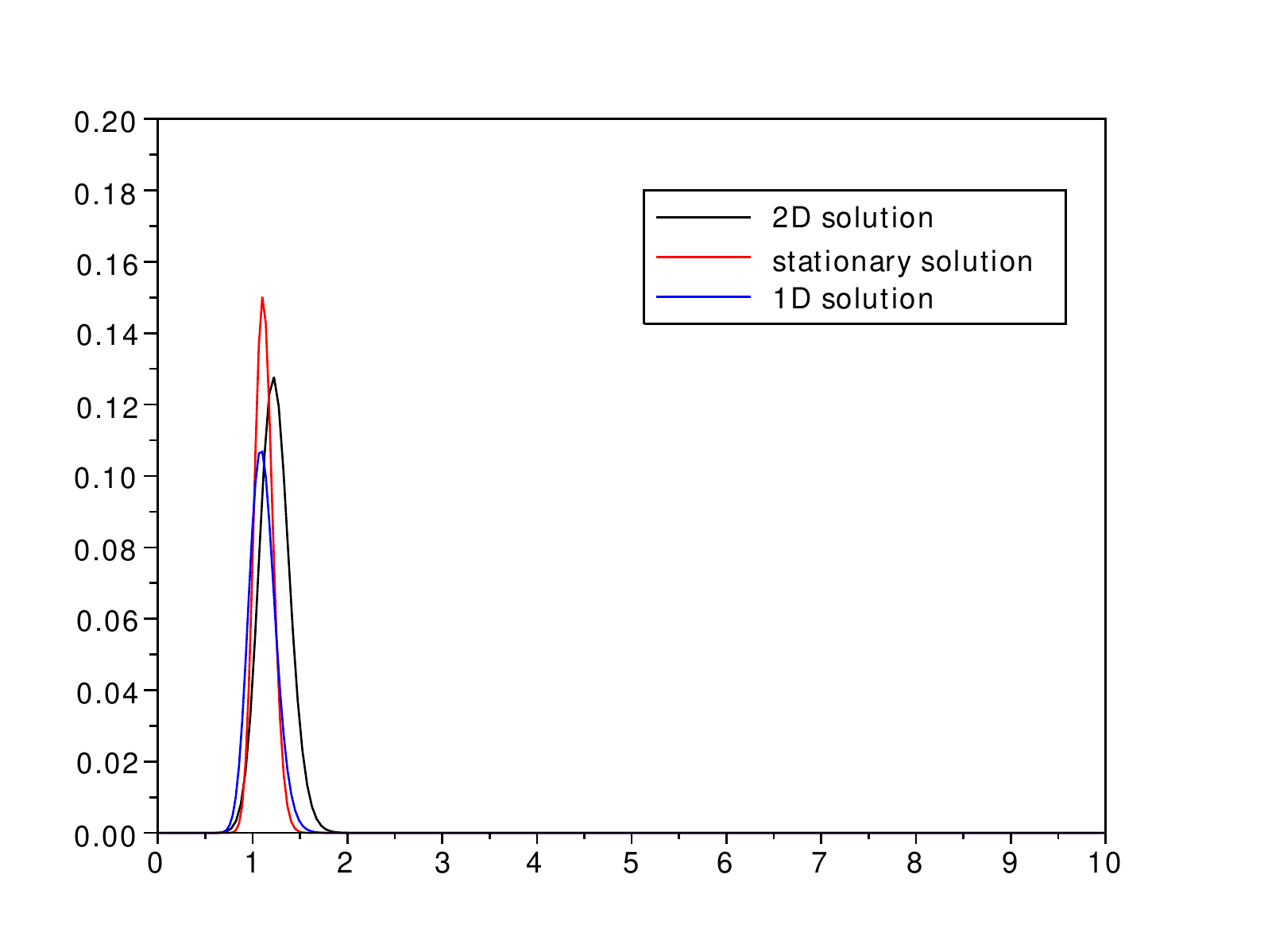}
\end{center}
\caption{Comparison of the solutions along $y$ for equations \eqref{FP} and \eqref{FP1}, and  with the stable state \eqref{SS}. Left~: unbiased. Right~: biased.}\label{marg}
\end{figure}

We can define the coefficient $\varepsilon$ as the ratio of the two eigenvalues associated to the Jacobian matrix of $F$, in such a way that $\varepsilon\ll 1$. This coefficient represents then the time scaling between the fast and slow variables. It is then possible to write the deterministic part of \eqref{EDS} as follows:
$$\begin{cases}
\varepsilon \textrm{d} x = f(x,y) \textrm{dt},\\
\textrm{d} y = g(x,y) \textrm{dt} .
\end{cases}$$
 Considering the limit of $\varepsilon$ going to zero, we can implicitly solve the first equation and define a curve $x^*(y)$ such that $f(x^*(y),y)=0$. Replacing in the equation for the slow variable $y$, we get:
$$\dot y =g(x^*(y),y).$$  
Considering now the stochastic term, we end up with the stochastic differential equation:
\begin{equation}\label{EDS1}
\dot y = g(x^*(y),y)+ \beta d\xi .
\end{equation}
We can then consider the associated  partial differential equation on $[0,\infty]\times \Omega_y$, with $\Omega_y=[-y_{m},y_{m}]$:
\begin{equation}\label{FP1}
\partial_t q + \partial_y \left( g(x^*(y),y) q - \frac{\beta^2}{2} \partial_y q\right)=0,
\end{equation}
where $q=q(t,y)$ is the probability distribution function representing the probability that at time $t\geq 0$, the firing rate is in $y$. This is a one-dimensional Fokker-Planck equation, and we can endowed it by means of no flux boundary conditions, for  $y=\{ -y_{m},y_{m}\}$:
$$g(y) q - \frac{\beta^2}{2} \partial_y q = 0, $$
and the normalized initial condition $q(0,y)=q_0(y)$, which is the projection of $p(t,\nu_1,\nu_2)$ only along the $y$ variable. The slow behavior of the solution persists, since $q$ lives on the equilibrium manifold along the slow direction. Nevertheless, computational time costs are reduced by using implicit in time numerical schemes. Moreover, we can compute an approximation of the potential function $V$ and of the stable state. In fact, for a one-dimensional Fokker-Planck equation the stable state is given by:
\begin{equation}\label{SS}
q_s(y)=\exp \left( \frac{-2G(y)}{\beta^2}\right),
\end{equation}
where $G(y)$ is the potential function associated to $g(x^*(y),y)$ and defined by:
$$G(y)=-\int g(x^*(z),z)\ dz.$$
As shown in figure \ref{marg}, the complexity reduced equation \eqref{FP1} of the initial Fokker-Planck model \eqref{FP} gives very good results, both in the unbiased and biased cases.
Therefore, it is possible to compute the wanted macroscopic quantities: reaction times and performance, as done in \citet{CCDM}. 

\section{Application to a three-well potential}
\begin{figure}[H]
\begin{center}
\includegraphics[width=.49\textwidth]{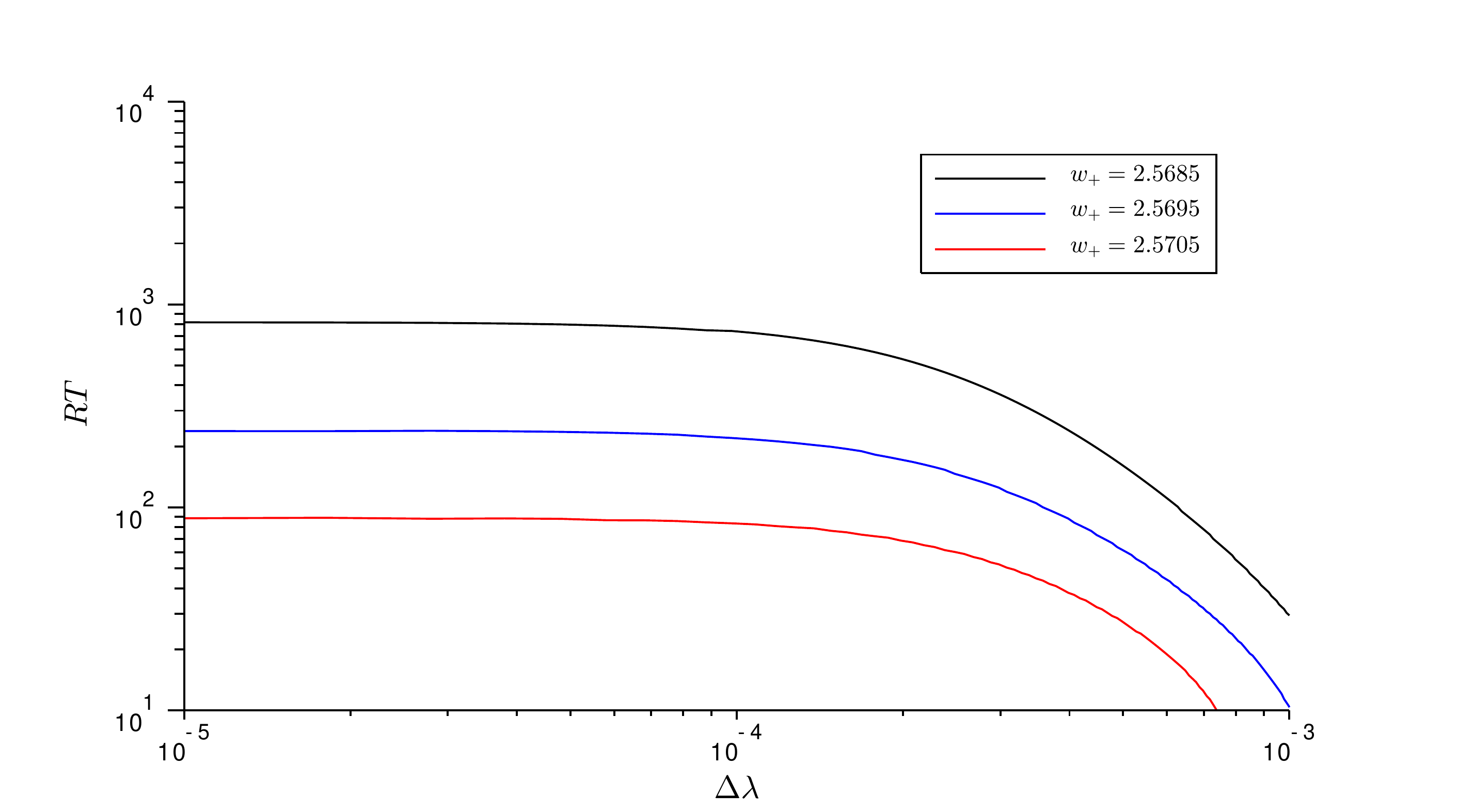}\; 
\includegraphics[width=.49\textwidth]{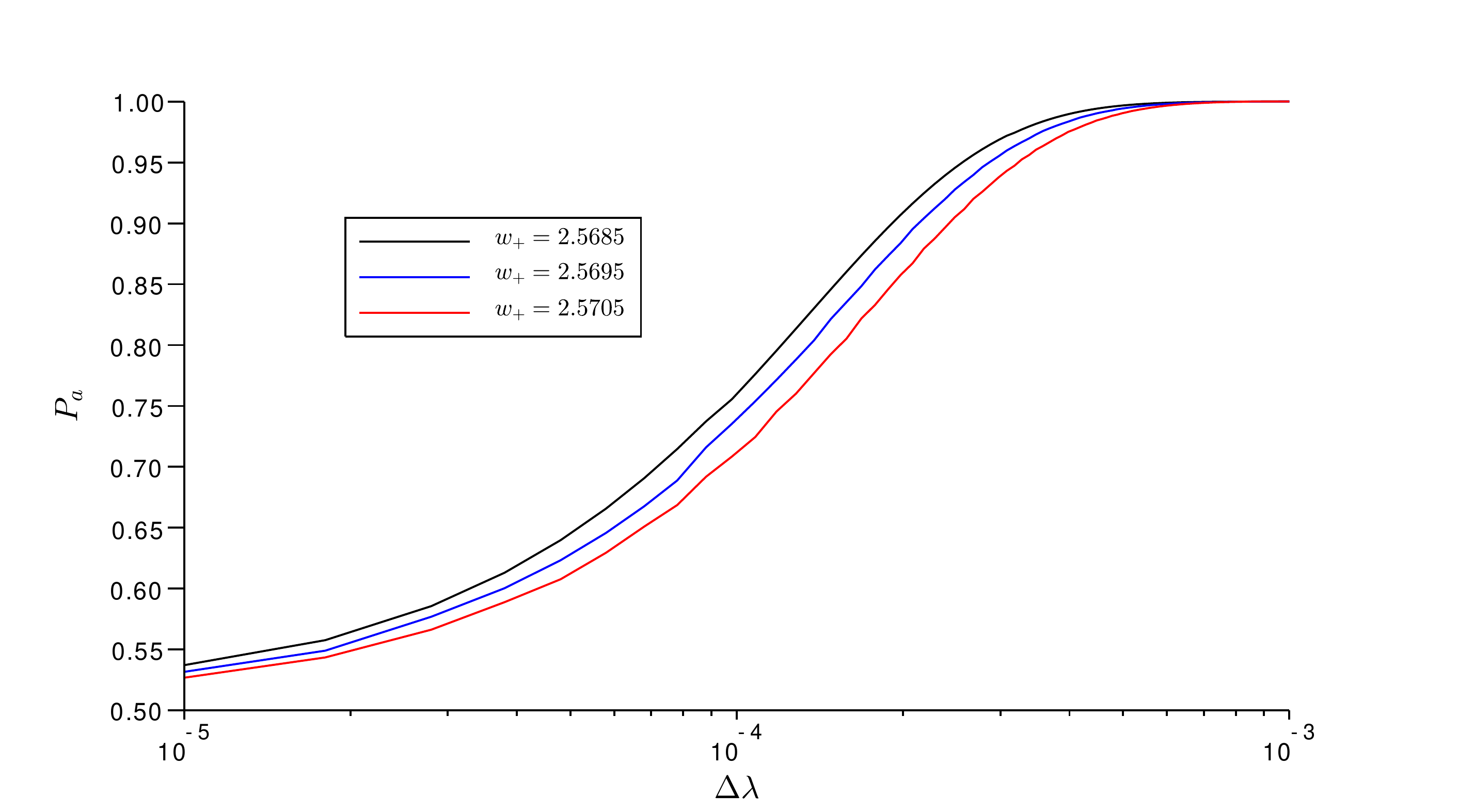}
\end{center}
\caption{Left~: Reaction times. Right~: Performance.}\label{rt-perf}
\end{figure}

The complexity reduction we have discussed for the double well potential is also useful for studying more complex situations like a three well potential. The main difference with what has been done previously is in a modified definition of the response function $\phi$. This situation is more realistic of what happens in visual decision making. Before a decision is made neurons firing rates are all concentrated at around a certain frequency value (the middle well) and they do migrate towards the other values (external wells) when the decision is made. 
Reaction times in the three well potential case are given by the exit times from the middle well to get to one of the external wells (the deeper one in the biased case). Whereas performance is defined by the density being, at equilibrium or for large times, into one specific well. 
In figure \ref{rt-perf} we plot the computed reaction times (left) and performance (right) both with respect to the difference on the applied stimuli $\Delta \lambda$ and for different values of the  coefficient $w_+$, which is one of connectivity coefficients used in the definition of the weights $w$ and $\hat{w}$.
 
The chosen values of the coefficient $w_+$ corresponds to a sub-critical bifurcation situation: the system passes form three minima to two minima. 
Indeed, for $w_+=2.5685$ the underlying potential has three stable equilibrium points and two unstable ones, for  $w_+=2.5695$ the middle well gets flat (the system is at the bifurcation) and for 
$w_+=2.5705$ the underlying potential has two wells and separated by a maximum which has replaced the middle well. Concerning reaction times (left), the larger the bias is the faster firing rates move towards the deeper well, and when the middle well disappears it become easier to take take a decision. Concerning, performance (right), the larger is the bias the more the subjects will give the expected answer, but the disappearance of the middle point doesn't increase the density of  the good decision since the bias also implies higher maximum values in the potential to overcome in order to get to the expected well. The same behavior was obtained in \citet{RL} for a similar problem and in several experimental results. Nevertheless, the approach proposed in \citet{RL} is valid only in a neighborhood of the spontaneous state (middle well), whereas the analysis and results presented here are valid on the whole domain of definition of the problem. 

\section{Conclusions}
We have discussed here in the framework of bi-stability view problems and of computational neurosciences, how the study of the partial differential equation associated to stochastic differential system of equation, can give complementary informations and can lead to the computation of macroscopic quantities (as reaction times and performance) of interest in the modeling of interacting population of neurons. In particular, we have presented the complexity reduction method based on the slow-fast behavior of the solution of a given stochastic differential system and applied it to a three well potential case.



\begin{thebibliography}{99}

\small 

\bibitem[Berglund $\&$ Gentz(2005)]{BG} {\sl N. Berglund and B. Gentz},
 Noise-induced phenomena in slow-fast dynamical systems: a sample-paths
approach. In: Probability and its applications. Springer, New York, (2005).

\bibitem[Carrillo, Cordier $\&$ Mancini (2011)]{CCM} 
{\sl J. A. Carrillo, S. Cordier and S. Mancini}, 
A decision-making Fokker-Planck model in computational
neuroscience, J. Math. Biol., 63, pp. 801-830, (2011).\\
doi: 10.1007/s00285-010-0391-3 

\bibitem[Carrillo, Cordier $\&$ Mancini (2013)]{CCM2} 
{\sl J. A. Carrillo, S. Cordier and S. Mancini},
One dimensional Fokker-Planck reduced dynamics of decision making models in Computational Neuroscience, Commun. Math. Sci., 11(2), pp. 523-540, (2013).\\
doi: 10.4310/CMS.2013.v11.n2.a10

\bibitem[Carrillo, Cordier, Deco $\&$ Mancini (2013)]{CCDM} 
{\sl J. A. Carrillo, S. Cordier, G. Deco and S. Mancini},
General One-Dimensional Fokker-Planck Reduction of Rate-equations
models for two-choice decision making, PLoS ONE 8(12): e80820. (2013)\\
doi:10.1371/journal.pone.0080820.

\bibitem[Deco $\&$ Martì (2007)]{DM} 
{\sl G. Deco and D. Martì},
Deterministic Analysis of Stochastic Bifurcations in Multi-stable Neurodynamical Systems, Biol. Cybern., 96(5), pp. 487-496, (2007).\\
doi: 10.1007/s00422-007-0144-6

\bibitem[Galam(2012)]{Galam} 
{\sl S. Galam}, Sociophysics : A Physicist's Modeling of Psycho-political Phenomena. In: Understanding Complex Systems.  Springer-Verlag: Berlin, (2012).

\bibitem[Roxin $\&$ Ledberg (2008)]{RL}
{\sl A. Roxin and A. Ledberg},
Neurobiological models of two-choice decision making can be reduced to a one-dimensional nonlinear diffusion equation, PLoS Comput Biol 4(3): e1000046. (2008).\\
doi:10.1371/journal.pcbi.1000046

\bibitem[Wilson $\&$ Cowan (1972)]{WC} 
{\sl H. R. Wilson and J. D. Cowan},
Excitatory and inhibitory interactions in localized populations of model neurons,
Biophy. J., 12(1), pp. 1-24, (1972).

\end{thebibliography}

\end{document}